\documentclass[10pt,a4paper]{amsart}

\usepackage{amsmath,amsthm,amstext,amsfonts,float,amssymb,mathrsfs}
\usepackage{mathptmx}
\usepackage[all]{xy}

\title[An application of Guillemin-Abreu theory]
{An application of Guillemin-Abreu theory \\ to a non-abelian group
action}

\author{Aleksis Raza}

\address{Department of Mathematics, 180 Queens Gate, Imperial College London SW7 2BT UK}

\email{aleksis.raza@imperial.ac.uk}



\newtheorem{defin}{Definition}

\newtheorem{theorem}[defin]{Theorem}
\newtheorem{prop}[defin]{Proposition}
\newtheorem{coro}[defin]{Corollary}

\newcommand{\bbR}{\mathbb{R}}
\newcommand{\bbC}{\mathbb{C}}
\newcommand{\tbbC}{\mathbb{C}}

\newcommand{\bbP}{\mathbb{P}}

\newcommand{\bbZ}{\mathbb{Z}}
\newcommand{\zi}{{\rm i}}

\newcommand{\wC}{\widehat{\mathbb{C}}}

\newcommand{\delbar}{{\overline{\partial}}}

\newcommand{\U}{{\rm U}}
\newcommand{\SU}{{\rm SU}}

\renewcommand{\ln}{\log}
\renewcommand{\geq}{\geqslant}

\begin{document}

\maketitle

\begin{abstract}
This note is a step towards demonstrating the benefits of a
symplectic approach to studying equivariant K\"ahler geometry. We
apply a local differential geometric framework from K\"ahler toric
geometry due to Guillemin \& Abreu to the case of the standard
linear $\SU(n)$ action on $\bbC^n\setminus\{0\}$. Using this
framework we (re)construct a scalar-flat K\"ahler metric on the
blow-up of $\bbC^n$ at the origin from data on the moment polytope.
\end{abstract}

\section{Introduction}

A symplectic toric manifold $(M,\omega,\tau,\mu)$ is a compact
symplectic manifold of dimension $2n$ equipped with an effective
hamiltonian action $\tau$ of the real $n$-torus
$T^n=\bbR^n/2\pi\bbZ^n$ and moment map $\mu:M\to(\bbR^n)^*$. By
Delzant theory $(M,\omega,\tau,\mu)$ admits a canonical
$\omega$-compatible $T^n$-invariant complex structure $J$,
\cite{del88:ham_per_imag_con_app_mom}. Hence every symplectic toric
manifold is canonically a K\"ahler toric manifold. Furthermore,
symplectic toric manifolds are completely classified by Delzant
polytopes $\Delta\subset (\bbR^n)^*$ i.e.\ moment polytopes
satisfying certain additional integrality conditions.

According to a differential geometric construction in toric
geometry one can in fact encode all the $T^n$-invariant K\"ahler
geometry of a symplectic toric manifold in terms of data on its
Delzant moment polytope $\Delta$. This construction relies on the
interplay between complex and symplectic structures in K\"ahler
geometry in the following sense. The data we refer to is a family
of certain smooth functions on $\Delta$ that determine every
possible $\omega$-compatible $T^n$-invariant complex structure on
$M$. These functions are obtained via a Legendre coordinate
transform from complex (holomorphic) coordinates on $M$ to
symplectic (action-angle) coordinates on $\Delta$. This coordinate
transform identifies K\"ahler potentials $f$ over $M$ (which
determine $T^n$-invariant $J$-compatible symplectic structures on
$M$ within a fixed cohomology class) as Legendre duals to {\it
symplectic potentials} $g$ on $\Delta$ (which determine
$T^n$-invariant $\omega$-compatible complex structures on $M$
within a fixed diffeomorphism class),
\cite{guil94:kah_str_tor_var, abr03:kah_geom_tor_man_sym_coor}. We
refer to this construction as Guillemin-Abreu theory.

In this paper we show how Guillemin-Abreu theory can be applied to
the case of the standard linear $\SU(n)$-action on $\bbC^n \setminus
\{0\}$. Using this construction we describe an interesting example.
We (re)construct a $\U(n)$-invariant, scalar-flat, K\"ahler metric
on $\wC^n$, the blow-up of $\bbC^n$ at the origin. This metric was
originally identified by Simanca who generalized the well-known
Burns metric on $\wC^2$, \cite{sim91:kah_met_con_scal_cur_bun_CP}.
The main purpose of this note is to illustrate that, in spirit of
Guillemin \& Abreu's work, doing K\"ahler geometry in symplectic
coordinates as opposed to the usual complex coordinates makes the
formulae quite elegant and the calculations more manageable.
Consequently the symplectic setting might be more appropriate for
working with these metrics, as in \cite{AP04:blow_des_kah_orb_csck}
which uses the Burns-Simanca metric to construct constant scalar
curvature K\"ahler metrics on blow-ups, for example.

I am grateful to my thesis examiners Dr. Andrew Dancer and Dr.
Richard Thomas for their comments and to my thesis supervisor
Prof. Simon Donaldson for his guidance.

\section{Guillemin-Abreu theory} \label{gat}

Let $(M,\omega,J,\tau,\mu)$ be a K\"ahler toric manifold and
$\Delta=\mu(M)$ its Delzant polytope. We first describe the local
structure of $M$. Consider the dense, open subset $M^\circ \subset
M$ where the $T^n$-action is free. Let $T^n_\bbC=\bbC^n/2\pi \zi
\bbZ^n =\bbR^n \times \zi T^n=\{w=a+\zi b: a \in \bbR^n, b \in
T^n\}$. $(M^\circ,J)\cong T^n_\bbC$ i.e. $(a,b)$ are complex
(holomorphic) coordinates on $M^\circ$ (see Appendix A of
\cite{abr03:kah_geom_tor_man_sym_coor}). The $T^n$-action on
$M^\circ$ is $(t,w)\mapsto w+\zi t$. The K\"ahler form $\omega$ is
given by $2\zi\partial \delbar f$ where $f\in C^\infty(M^\circ)$.
Since $\omega$ is $T^n$-invariant $f=f(a)\in C^{\infty}(\bbR^n)$.
Thus
\begin{equation} \label{om}
\omega=\sum_{i,j=1}^n \frac{\partial^2 f}{\partial a_j \partial
a_k}da_i\wedge db_j.
\end{equation}
We now describe the interior $\Delta^\circ$ of $\Delta$. Suppose
$\Delta $ consists of $d$ facets (codimension-1 faces). Then $
\Delta=\{x \in (\bbR^n)^*: \langle x,u_i \rangle \geq \lambda_i,
\; i=1,\dots,d\}$ where $u_i$ is the integral primitive inward
pointing normal vector to the $i$th facet of $\Delta$ (see
\cite{guil94:kah_str_tor_var}). Let $l_i$ be affine functions on
$(\bbR^n)^*$ defined by $ l_i: x \mapsto \langle x,u_i \rangle -
\lambda_i$. Then $ x \in \Delta^{\circ} $ if and only if $ l_i(x)
> 0 $.
\begin{theorem}[\cite{guil94:kah_str_tor_var}] \label{g1}
$\mu$ factors into
\[
\xymatrix{ \bbR^n \times \zi T^n \ar[r]^\mu \ar[d]_p & (\bbR^n)^* \\
\bbR^n \ar[ur]_{\mu_f}}
\]
and this diagram commutes. Moreover, $\mu_f$ is the Legendre
transform $a \mapsto df_a=x$ associated to $f$ and is a
diffeomorphism onto $\Delta^\circ\subset (\bbR^n)^*$. Furthermore,
there exists an inverse Legendre transform
\[
\begin{array}{cccc}
\mu_f^{-1}: & \Delta^\circ & \to & \bbR^n, \\
  & x & \mapsto & dg_x=a
\end{array}
\]
where the function $g \in C^\infty(\Delta^\circ)$ is the Legendre
dual to $f$ i.e.
\begin{equation}\label{ldual} f(a) + g(x) = \sum_{i=1}^n
\frac{\partial f}{\partial a_i}\frac{\partial g}{\partial x_i}.
\end{equation}
\end{theorem}
Guillemin's set-up provides us with a means of encoding the K\"ahler
data on $M^\circ$, originally given in terms of $f$ and the
coordinates $(a,b)$, into symplectic (action-angle) coordinates
$(x,y)$ on $\Delta^\circ \times T^n$ and $g$ through the map $(a,b)
\mapsto (x,y)$ which is the Legendre transform $\mu_f$ on the first
factor and the identity on the second factor. Guillemin introduces
the function $g(x)=\frac{1}{2}\sum_{i=1}^d l_i(x)\ln l_i(x)$. The
particular form of this function guarantees that it is convex and
smooth on $\Delta^\circ$ and has the appropriate singular behavior
on the boundary $\partial \Delta$ of $\Delta$. It is the Legendre
dual of the K\"ahler potential $f$ that defines the canonical
$T^n$-invariant K\"ahler metric $\omega(\cdot,J\cdot)$ on $M^\circ$.
Using Guillemin's set-up one can construct any other $T^n$-invariant
K\"ahler metric on $M$ in the class $[\omega]$ purely from the
combinatorial data on $\Delta$ employing such functions $g \in
C^\infty(\Delta^\circ)$. Those $g$ whose Legendre duals $f$ define
$T^n$-invariant K\"ahler metrics can be regarded as `potentials' for
complex structures on $M$ in analogy to the K\"ahler potentials $f$
for symplectic structures on $M$.

Just as symplectic structures within a fixed cohomology class are
parameterized in terms of K\"ahler potentials $f$ through the
$\partial\overline{\partial}$-lemma, Abreu provides an analogous
theorem for parameterizing complex structures within a fixed
diffeomorphism class in terms of symplectic potentials $g$:
\begin{theorem}[{\rm \cite{abr03:kah_geom_tor_man_sym_coor}}] \label{a1}
Let $(M_\Delta,\omega_\Delta)$ be a symplectic toric $2n$-manifold
corresponding to the Delzant polytope $\Delta$ which has $d$
facets. Then a $T^n$-invariant, $\omega_\Delta$-compatible complex
structure $J$ on $M_\Delta$, given at a point in the coordinates
$(x,y)$ by
\begin{equation}\label{ggi}
\begin{pmatrix}
0 & -G^{-1} \\
G & 0
\end{pmatrix},
\end{equation}
is determined by a smooth function
\begin{equation}\label{ggg}
g(x)=\frac{1}{2} \sum_{i=1}^d l_i(x) \log l_i(x) + h(x)
\end{equation}
on $\Delta^\circ$, where $ h(x) \in C^\infty(\Delta)$ (i.e. there is
an open set $U \subset (\bbR^n)^*$ containing $\Delta$ and an
$\widetilde{h}\in C^\infty(U)$ which restricts to $h$ on $\Delta$),
the hessian matrix $G$ of \eqref{ggg} is positive definite on
$\Delta^\circ$ and
\begin{equation}\label{delta}
\det G^{-1}=\delta(x)\prod_{i=1}^d l_i(x)
\end{equation}
where $\delta(x) \in C^\infty(\Delta)$ and is strictly positive on
$\Delta$. Conversely, every $g$ of the form \eqref{ggg} determines a
$T^n$-invariant, $\omega_\Delta$-compatible complex structure on
$(M_\Delta,\omega_\Delta)$ which in the $(x,y)$ coordinates is of
the form \eqref{ggi}.
\end{theorem}
The Guillemin-Abreu paradigm is that one can recover $f$ from $ g$
so it is enough to work with $g$ and data on ${\Delta}$. One final
result to recall from the Guillemin-Abreu framework is that Abreu
applies the above ideas to formulate an elegant expression for the
scalar curvature of the K\"ahler metric defined by \eqref{om}.
Furthermore, Abreu also provides `symplectic' extremal K\"ahler
condition. These are summed up in the following
\begin{theorem}[\cite{abr98:kah_geom_tor_var_ext_met}]
The scalar curvature of \eqref{om} is
\begin{equation}\label{scalt}
S(g) = - \frac{1}{2}\sum^{n}_{i,j=1} \frac{\partial^2
G^{ij}}{\partial x_i \partial x_j}
\end{equation}
where $G^{ij}$ is the $(i,j)$th entry of $G^{-1}$. Furthermore,
this K\"ahler metric is extremal if and only if
\begin{equation}\label{excon}
\frac{\partial S}{\partial x_i}={\rm constant}
\end{equation}
for $i=1,\dots,n$ i.e. $S$ is an affine function of $x$.
\end{theorem}

\section{A non-abelian group action} \label{two}

Let $(M,\omega,J)$ be a K\"ahler $n$-fold. The scalar curvature
$S_J$ of the K\"ahler metric $\omega(\cdot,J\cdot)$ is given by
\begin{equation} \label{eqsc}
S_J \omega^n=n!\Theta_J \wedge \omega^{n-1}
\end{equation}
where $\Theta_J$ is the Ricci form with respect to the complex
structure $J$. Consider $\tbbC^2\setminus \{0\}$ with standard
complex coordinates $z=(z_1,z_2)$ equipped with the standard
linear $\SU(2)$-action. K\"ahler potentials of $\SU(2)$-invariant
K\"ahler metrics are smooth functions on $\bbC^2\setminus\{0\}$ of
the form $f=f(s)$ where $s=|z_1|^2+|z_2|^2$ is the square of the
radius of the $\SU(2)$-orbits.  We are interested in studying
$\SU(2)$-invariant K\"ahler metrics on $\tbbC^2\setminus \{0\}$.
The standard approach is to fix the standard $J_0$ on
$\tbbC^2\setminus \{0\}$ and vary the symplectic structure using
the K\"ahler potentials $f(s)$ (through the `$\partial
\overline{\partial}$-lemma'). The scalar curvature of the K\"ahler
metric determined by an $f(s)$ i.e.
\begin{equation} \label{su}
h_{i\overline{j}} = \left[\frac{\partial^2 f}{\partial z_i
\partial \overline{z}_j}\right]_{i,j=1}^2 =
\begin{pmatrix}
f' + z_1\overline{z}_1 f'' &
z_2\overline{z}_1 f'' \\
z_1\overline{z}_2 f'' & f' +
z_2\overline{z}_2 f'' \\
\end{pmatrix},
\end{equation}
is then deduced from \eqref{eqsc}. Not all $f(s)$ determine
$\SU(2)$-invariant K\"ahler metrics on $\bbC^2\setminus \{0\}$.
$h_{i\overline{j}}$ is positive definite if and only if $ f'(s)>0,
f''(s) > -s^{-1}f'(s)$. The scalar curvature of \eqref{su}
computes to be a complicated expression in $f$, see
\cite{cal82:ext_kah_met, sim91:kah_met_con_scal_cur_bun_CP}.

\subsection{Twisted Guillemin-Abreu theory}

The key observation is that the $\SU(2)$-invariant K\"ahler
metrics on $\tbbC^2\setminus \{0\}$ are also invariant under
$\U(2)$ and in particular under $T^2\cong \U(1) \times \U(1)
\subset \U(2)$. We conclude from this observation that it is
viable to employ Guillemin-Abreu theory for this standard linear
$\SU(2)$-action. Set $ w_j = \ln z_j $ where $ w_j=a_j+\zi b_j$, $
j=1,2 $ so $s= |z_1|^2+|z_2|^2 = e^{w_1}e^{\overline{w}_1} +
e^{w_2} e^{\overline{w}_2}=e^{2a_1}+e^{2a_2}$ parameterizes the
$\SU(2)$ orbits. Now we invoke the Guillemin-Abreu theory
introduced in the previous section. Applying the Legendre
transform associated to a choice of $f(s)$ gives
\[
x_i=\frac{\partial f}{\partial a_i}= f'(s) \frac{\partial
s}{\partial a_i}=2 e^{2a_i} f'(s),
\]
$i=1,2$ and $ x_1+x_2=2(e^{2a_1}+e^{2 a_2})f'=2sf'=\gamma(s)$. Let
$h$ be the inverse function to $\gamma$, then $ s= h(x_1+x_2)$.
Substituting $ e^{2 a_1}=e^{2 a_2} x_1x_2^{-1}$ into
$s=e^{2a_1}+e^{2a_2}$ gives $ s=e^{2a_2} + e^{2a_2}
x_1x_2^{-1}=e^{2 a_2}(x_1+x_2)x_1^{-1}$ i.e. $ e^{2 a_1}=s
x_1(x_1+x_2)^{-1}$. Similarly $ e^{2a_2}=s x_2(x_1+x_2)^{-1}$. The
Legendre dual of $f$ is given by \eqref{ldual} in Theorem
\ref{g1}. We have established above that $s=h(x_1+x_2)$ hence
$f=f(s)=f(h(x_1+x_2))=f(x_1+x_2)$. A brief calculation shows that
\begin{equation} \label{ssp}
g(x)= \frac{1}{2} \left( x_1 \ln x_1 + x_2 \ln x_2 + F(x_1+x_2)
\right)
\end{equation}
with
\begin{equation}\label{FF}
F(x_1+x_2)=F(t)= t \ln \left(h(t)t^{-1}\right)- 2f(s(t))
\end{equation}
where we have set $t=x_1+x_2$. The Guillemin-Abreu approach is to
fix the standard $\omega_0$ on $\tbbC^2\setminus \{0\}$ and vary
the complex structure using the symplectic potentials $g(x)$
(through Theorem \ref{a1}). The standard $\omega_0$ is given by
$f(s)=s/2$. Hence the moment polytope $\Delta_{\bbC^2}$ for the
standard $T^2$-action is the positive orthant $\bbR^2_{\geqslant
0} \subset \bbR^2$ with standard symplectic (action) coordinates
$(x_1,x_2)=(|z_1|^2,|z_2|^2)$. $g(x)$ is a smooth function on
$\Delta_{\bbC^2}$ while $F(t)$ is a smooth function on
$(0,\infty)$. We call $F(t)$ the {\it $t$-potential} of $g$ since
it is the `$t$-part' of the symplectic potential $g$. The hessian
matrix of $g$ is
\[
G=(G_{ij})=\frac{1}{2} \left(\begin{array}{cc}\frac{1}{x_1}+F''(t)
& F''(t) \\ F''(t) & \frac{1}{x_2}+F''(t)
\end{array}\right)
\]
and $ \det G = 4^{-1}\left((x_1x_2)^{-1}+ t(x_1x_2)^{-1}
F''\right)$. It follows from the discussion in the previous
section that $G$ must be positive definite. This in turn implies
that $F''(t)>-t^{-1} $ and, furthermore, this condition is also
sufficient i.e. only functions of the form $F(t)$ that satisfy
this property determine $\SU(2)$-invariant K\"ahler metrics on
$\bbC^2\setminus \{0\}$. To summarize
\begin{prop} \label{propsuinv}
The K\"ahler metric \eqref{su} has a symplectic potential given by
\eqref{ssp}. Conversely, any such smooth function on $\bbR^2_{\geq
0}$ with $F$ satisfying $F''(t)>-t^{-1} $ determines such a
K\"ahler metric.
\end{prop}
The inverted hessian of $g$ is
\begin{equation}\label{Gi}
G^{-1}=(G^{ij})=\frac{2x_1x_2}{1+ tF''}\left(\begin{array}{cc}
\frac{1}{x_2}+F''(t)
& - F''(t) \\
-F''(t) & \frac{1}{x_1}+F''(t)
\end{array}\right).
\end{equation}
Applying Abreu's scalar curvature formula \eqref{scalt} to
\eqref{Gi} shows (see \cite{Raz04:scal_cur_mul_free_act}) that
\begin{prop} \label{propode}
The scalar curvature of the K\"ahler metric \eqref{su} is given by
\[
S(g)=t^{-1} \left(t^3F''(1 + tF'')^{-1}\right)''.
\]
\end{prop}
Propositions \ref{propsuinv} and \ref{propode} hold true for
arbitrary dimension. Let $z=(z_1,\dots,z_n)$ be standard complex
coordinates on $\bbC^n\setminus \{0\}$. Let
$s=\sum_{i=1}^n|z_i|^2$. An $\SU(n)$-invariant K\"ahler metric on
$\bbC^n \setminus \{0\}$ is given by
\begin{equation} \label{sun}
h_{i\overline{j}} = \left[\frac{\partial^2 f}{\partial z_i
\partial \overline{z}_j}\right]_{i,j=1}^n =
\left[f'\delta_{ij}  + z_i\overline{z}_jf''\right]_{i,j=1}^n.
\end{equation}
Repeating the construction we carried out for the case $n=2$ we
obtain symplectic coordinates $x = (x_1,\dots,x_n)$ on the
positive orthant $\bbR^n_{\geq 0}=\Delta_{\bbC^n}$ i.e. the moment
polytope for the standard $T^n$-action on $\bbC^n$ (with its
standard symplectic structure). One simply has to work through the
algebra and extend the identities already derived earlier to the
$n$ case. Thus
\begin{theorem}[\cite{Raz04:scal_cur_mul_free_act}] \label{myfirst}
The K\"ahler metric \eqref{sun} has a symplectic potential given
by
\begin{equation}  \label{spp}
g(x)=\frac{1}{2} \left[\sum_{i=1}^n x_i \ln x_i + F(t)\right]
\end{equation}
with
\begin{equation} \label{F}
F(t)= t \ln \left(s(t)t^{-1}\right) - 2f(s(t))
\end{equation}
where $ t=\sum_{i=1}^nx_i=2sf'(s)$ and its scalar curvature is
given by
\begin{equation} \label{scn}
S(g)=t^{1-n} \left(t^{n+1}F''\left(1+tF''\right)^{-1}\right)''.
\end{equation}
Conversely, any function of the form \eqref{spp} on
$\Delta_{\bbC^n}$ with \eqref{F} satisfying $F''(t)>-t^{-1}$
determines such a K\"ahler metric.
\end{theorem}

\section{Applications of Theorem \ref{myfirst}}

As a first application of Theorem \ref{myfirst} we verify the
well-known result that the scalar curvature of the Fubini-Study
metric on $\bbC\bbP^n$ is constant. The $t$-potential of the
Fubini-Study metric on $ \bbC\bbP^n$ is $
F_{\bbC\bbP^n}(t)=(1-t)\ln(1-t)$. Substituting $F_{\bbC\bbP^n}(t)$
into \eqref{scn} shows that
\begin{coro}
The Fubini-Study metric on $\bbC\bbP^n$ has scalar curvature
$n(n+1)$.
\end{coro}
Let $\wC^n$ denote the blow-up of $\bbC^n$ at the origin. Recall
that the Burns metric on $\widehat{\bbC}^2$ is the restriction of
the standard product metric on the ambient space $\bbC^2 \times
\bbC\bbP^1$ when $\wC^2$ is considered as a hypersurface in
$\bbC^2 \times \bbC\bbP^1$
\cite{sim91:kah_met_con_scal_cur_bun_CP}. We refer to the
restriction of the standard product metric on $\bbC^n \times
\bbC\bbP^{n-1}$ to $\widehat{\bbC}^n$ as the {\it generalized
Burns metric}. In light of this we are led to ask whether the
generalized Burns metric on $\widehat{\bbC}^n$ is also
scalar-flat. Unsurprisingly this is not the case. The
$t$-potential of the generalized Burns metric is $
F_{\widehat{\bbC}^n}(t)=(t-1) \ln (t-1) -t \ln t - t+1$.
Substituting this into \eqref{scn} shows that the scalar curvature
of the generalized Burns metric on $\widehat{\bbC}^n$ is
$S(g_{\widehat{\bbC}^n})=(n^2-3n+2)t^{-2}$. So when do
$t$-potentials of this form give rise to $\SU(n)$-invariant
scalar-flat K\"ahler metrics of this form? Setting
$S(g_{\widehat{\bbC}^n})$ to zero gives the quadratic $n^2-3n+2=0$
whose solutions are $n=1$ and $n=2$. Thus
\begin{coro} \label{nsf}
The generalized Burns metric on $\wC^n$ is scalar-flat in and only
in dimension $1$ and $2$ i.e. regarding $\widehat{\bbC}^n$ as a
hypersurface in $\bbC^n \times \bbC\bbP^{n-1}$, the restriction of
the standard ambient product metric on this space to
$\widehat{\bbC}^n$ is scalar-flat in and only in dimension 1 and
2.
\end{coro}

\subsection{A scalar-flat K\"ahler metric on $\pmb{\widehat{\bbC}^n}$}

Consider now $\U(n)$-invariant K\"ahler metrics of zero scalar
curvature on $\bbC^n$. Then Theorem \ref{myfirst} leads to the ODE
\eqref{scn}$=0$ which solves (for $F''$) to give
\begin{equation}\label{red}
F''(t)=\frac{At+B}{t(t^n-(At+B))}
\end{equation}
where $A$ and $B$ are constants. Simanca proved the existence of a
$\U(n)$-invariant scalar-flat complete K\"ahler metric on the
total space of the bundle $\mathscr{O}(-1) \to \bbC\bbP^{n-1}$
i.e. the blow-up $\wC^n$ of $\bbC^n$ at the origin,
\cite{sim91:kah_met_con_scal_cur_bun_CP}. We shall now describe
this metric in the symplectic coordinate setting. Since the
K\"ahler metric we seek is also $T^n$-invariant we employ
Guillemin-Abreu theory, in particular Theorem \ref{a1}, to find
our boundary condition. The Delzant moment polytope
$\Delta_{\wC^n}$ corresponding to $\wC^n$ is the positive orthant
$\bbR^n_{\geq 0}$ with the vertex $p=(0,\dots,0)$ replaced by the
$n$ vertices $p+x_i$, $i=1,\dots,n$. We refer to $\Delta_{\wC^n}$
as the moment polytope of the 1-symplectic blow-up of $\bbC^n$ at
the origin (see \cite{Raz04:scal_cur_mul_free_act}). As a result
$\Delta_{\wC^n}$ has $(n+1)$ facets. Let $ l_1(x),\; l_2(x),\;
\dots, l_{n+1}(x)$ be the affine functions corresponding to these
facets i.e. $l_i(x)=x_i $ and $ l_{n+1}=\sum_{i=1}^nx_i - 1=t-1$.
That is, each affine function determines a hyperplane in $\bbR^n$
and these hyperplanes together trace out the boundary of
$\Delta_{\wC^n}$. The interior of $\Delta_{\wC^n}$ is $
\Delta^\circ_{\wC^n}=\{x\in \bbR^n:\; l_i(x)> 0, \;
i=1,\dots,n+1\}$. Let $g_{\wC^n}(x)$ be the symplectic potential
of the metric we seek. By Theorem \ref{myfirst}
\begin{equation}\label{tobe}
g_{\wC^n}(x) = \frac{1}{2}\left[\sum_{i=1}^nx_i\ln x_i +
F_{\wC^n}(t)\right]
\end{equation}
where $F_{\wC^n}(t)$ is the $t$-potential of $g_{\wC^n}(x)$.
Furthermore, for $g_{\wC^n}(x)$ to determine a scalar-flat
K\"ahler metric $F_{\wC^n}(t)$ must be of the form \eqref{red}.
The determinant of the hessian $G_{\wC^n}$ of $g_{\wC^n}(x)$ is $
\det G_{\wC^n} = 2^{-n} \left(1+ tF''_{\wC^n}(t)\right)
\prod_{i=1}^n x^{-1}_i$ and so
\begin{equation}\label{GG}
\det G^{-1}_{\wC^n}=2^n\prod_{i=1}^n x_i\left(1+
tF''_{\wC^n}\right)^{-1} = 2^n \prod_{i=1}^n x_i
(t^n-(At+B))t^{-n}.
\end{equation}
By Theorem \ref{a1} $\det G^{-1}_{\wC^n}$ should be of the form
\eqref{delta} with $l_i$ as given above and $\delta(x) \in
C^\infty(\Delta_{\wC^n})$ (in the sense described in Theorem
\ref{a1}) and positive. Writing \eqref{GG} as $ \det G^{-1}_{\wC^n}
= \delta_{\wC^n}(x) \prod_{i=1}^{n+1} l_i(x)$ such that
\begin{equation} \label{del}
\delta_{\wC^n}(x) = \frac{2^n(t^n-(At+B))}{t^n(t-1)}
\end{equation} gives us the
appropriate form of $\det G^{-1}_{\wC^n}$. We now have to find the
correct $A,B$ in \eqref{del} so that it satisfies the condition of
Theorem \ref{a1}. In its current form $\delta_{\wC^n}(x)$ becomes
singular at the boundary i.e. when $t \to 1$. Also $2^nt^{-n}>0$ and
smooth. Therefore $A,B$ must be such that $(t^n-At-B) \approx (t-1)$
for $t=1+\epsilon$ for small $\epsilon>0$. Then $
(1+\epsilon)^n-A(1+\epsilon)-B = 1 + n\epsilon +
O(\epsilon^2)-A(1+\epsilon)\approx (1+\epsilon)-1$ which gives $ (1
- (A+B)) + (n-A)\epsilon = \epsilon$. Hence $n-A=1$ and $1-(A+B)=0$
and we get $A=n-1,\;B=2-n$. It follows that $
(t^n-(n-1)t-(2-n))(t-1)^{-1} = \sum_{i=1}^{n-1}t^i-(n-2)$ which is
clearly smooth and positive on the whole of $\Delta_{\wC^n}$. Hence
$ \delta_{\wC^n}(t) = 2^n t^{-n} \left(\sum_{i=1}^{n-1}
t^i-(n-2)\right)$. As a result
\begin{equation}\label{burnsrules}
F''_{\widehat{\bbC}^n}(t)= \frac{(n-1)t+2-n}{t(t^n-(n-1)t-2+n)}
\end{equation}
such that $n>1$. We need to make sure that
$F''_{\widehat{\bbC}^n}(t)$ is non-singular for $t>1$ i.e.
$t^n-(n-1)t-2+n \neq 0$ if $t>1$. We have $ t^n-(n-1)t-(2-n) =
(t-1)\left(\sum_{i=1}^{n-1}t^i-(n-2)\right)$ and this is zero if
$(t-1)=0$ or if $ \sum_{i=1}^{n-1}t^i-(n-2) = 0$. It is clear that
for $t\geq 1$ and $n>1$, $ \sum_{i=1}^{n-1}t^i
> n-2$. It follows that \eqref{burnsrules} is non-singular for all
$t>1$ and hence we have
\begin{coro}\label{second}
There exists a $\U(n)$-invariant, scalar-flat, K\"ahler metric on
$\widehat{\bbC}^n$ determined by the symplectic potential
\eqref{tobe} on $\Delta_{\widehat{\bbC}^n}$ where
$F_{\widehat{\bbC}^n}(t) - (t-1)\ln(t-1)$ is a smooth function on
$[1,\infty)$ such that $F_{\wC^n}(t)$ satisfies
\eqref{burnsrules}.
\end{coro}
We refer to the K\"ahler metric determined by the symplectic
potential \eqref{tobe} as the Burns-Simanca metric. The case $n=2$
is special since in that case the metric on $\wC^2$ determined by
\eqref{tobe} is just the restriction of the product metric on
$\bbC^2 \times \bbC\bbP^1$. The crucial point is that for $n > 2$
the Burns-Simanca metric on $\wC^n$ is not the restriction of the
standard product metric on the ambient space $\bbC^n\times
\bbC\bbP^{n-1}$.

\subsubsection{Behaviour of the Burns-Simanca metric away from
the blow-up}

Let $d=(\delta_{ij})$ be the standard flat euclidean metric on
$\bbC^n$. Then a K\"ahler metric $h$ on a non-compact K\"ahler
$n$-fold $M$ is called asymptotically euclidean (AE) with rate of
decay $r^{-n}$, where $r$ is the radius function on $\bbC^n$, if
$h$ approximates $d$ under an appropriate biholomorphism between
$M$ (minus a compact subset) and $\bbC^n$.
\begin{prop}
The Burns-Simanca metric on $\wC^n$ is AE with rate of decay
$r^{2-2n}$.
\end{prop}
\begin{proof}
Let $(y_1,\dots,y_n)$ be the toric (angle) coordinates. The flat
metric on $\bbC^n\setminus \{0\}$ is given in the symplectic
coordinates $(x,y)$ by
\begin{equation}\label{fmm}
d_{(x,y)}= \sum_{i=1}^n\frac{1}{2x_i}dx_i\otimes dx_i + 2x_i
dy_i\otimes dy_i.
\end{equation}
The $(i,j)$th entries of the hessian matrix of \eqref{tobe} and
its inverse are
\begin{equation}\label{gbh}
(G_{\wC^n})_{ij}=\frac{1}{2}\left\{\begin{array}{ll} x^{-1}_i
+F''_{\wC^n}, & i=j \\  F''_{\wC^n}, & i \neq j
\end{array}\right.
\end{equation}
and
\begin{equation}\label{gbhi}
(G_{\wC^n})^{ij}=\frac{2}{1+t F''_{\wC^n}}
\left\{\begin{array}{ll} \left(1+ F''_{\wC^n}\sum_{k=1,k \neq i}^n
x_k \right)x_i, & i=j \\
-F''_{\wC^n} x_ix_j, & i \neq j,
\end{array}\right.
\end{equation}
respectively. See \cite{Raz04:scal_cur_mul_free_act}. By the
Guillemin-Abreu construction the Burns-Simanca metric is given by
the $2n\times 2n$ matrix $ h_{BS}= {\rm diag}
(G_{\wC^n},G^{-1}_{\wC^n})$. \eqref{gbh} can be written as the sum
$ G_{\wC^n} = A +B$ where $ A={\rm
diag}((2x_1)^{-1},\dots,(2x_n)^{-1})$ and $
B=((n-1)t+2-n)(2t(t^n-(n-1)t-2+n))^{-1}P$ where $P$ is an $n\times
n$ matrix of 1's. This gives us the upper left block of $h_{BS}$.
For \eqref{gbhi} the coefficient term is $2t^{-n}(t^n-(n-1)t+n-2)$
hence the off-diagonal terms of $G^{-1}_{\wC^n}$ are $
(G_{\wC^n})^{ij}= - 2((n-1)t+2-n)x_ix_jt^{-n-1}$. Note that in
\eqref{gbhi} the term $ \left(1+ F''_{\wC^n}(t)\sum_{k=1,k \neq
i}^n x_k \right)x_i = (1+(t-x_i)F''_{\wC^n}(t))x_i$. Hence the
diagonal terms in $G^{-1}_{\wC^n}$ are $ (G_{\wC^n})^{ii} = (2(1+
(t-x_i)F''_{\wC^n}(t))x_i)(1+tF''_{\wC^n}(t))^{-1}$ which using
the coefficient term gives $(G_{\wC^n})^{ii} =
2x_i((t^n-(n-1)t+n-2) + ((n-1)t+2-n)(t-x_i)t^{-1})t^{-n}$. Thus we
have that the inverted hessian matrix of \eqref{tobe} splits into
$ G^{-1}_{\wC^n} = C  + D$ where $C=t^{-n}(t^n-(n-1)t+n-2){\rm
diag}(2x_1,\dots, 2x_n)$ and $ D = 2t^{-n-1}((n-1)t+2-n)Q$ such
that $Q={\rm diag}(tx_1,\dots,tx_n) - [x_ix_j]_{i,j=1}^n$. This
gives us the lower right block of $h_{BS}$. Therefore we have $
h_{BS}= {\rm diag}(A,C) + {\rm diag}(B,D)= E+F$ so that $
E=\sum_{i=1}^n\frac{1}{2x_i}dx_i\otimes dx_i +
2x_it^{-n}(t^n-(n-1)t+n-2)dy_i\otimes dy_i$. Now we use the flat
metric coordinates $(\lambda,\mu)$. Set
\[
\lambda_i=\sqrt{2x_i}\cos y_i,\; \mu_i= \sqrt{2x_i}\sin y_i,
\]
$i=1,\dots,n$. Then the flat metric \eqref{fmm} becomes of the
standard form $ d_{(\lambda,\mu)}=\sum_{i=1}^n d\lambda_i\otimes
d\lambda_i + d \mu_i \otimes d\mu_i$. The component $E$ of
$h_{BS}$ in these $(\lambda,\mu)$ coordinates is given by $ E=
\sum_{i=1}^n d\lambda_i\otimes d\lambda_i + d \mu_i \otimes d\mu_i
+ O(u^{-n+1})$ where $ u=2^{-1}\sum_{i=1}^n\lambda^2_i + \mu^2_i$
is the variable $t$ in the $(\lambda,\mu)$ coordinates.  For the
component $F$ of $h_{BS}$ its clear that for $u\to\infty$ the
coefficient terms are $O(u^{-n})$. It follows that as $t=u \to
\infty$, $h_{BS}=d_{(\lambda,\mu)}+O(u^{-n+1})$. Since by
construction the coordinate $u$ represents the square of the
radius function $r$ on $\bbC^n$ we deduce that $u^{1-n}\equiv
r^{2(1-n)}$.
\end{proof}

\bibliography{mybib}

\end{document}